\documentclass[12pt]{amsart}
\usepackage{latexsym}
\usepackage{psfig}
\usepackage{amsmath,amssymb}

\newtheorem{thm}{Theorem}[section]

\newtheorem{Lemma}[thm]{Lemma}
\newtheorem{Remark}[thm]{Remark}
\newtheorem{cl}[thm]{Claim}

\numberwithin{equation}{section}

\DeclareMathOperator{\m}{mod}

\author{Joseph Cohen}
\title{Primitive Roots in Quadratic Fields II}
\email{coheny@techunix.technion.ac.il}
 \address{Math Department,
Technion, Haifa, 32000, Israel}
\date{\today}
\thanks{In partial fulfillment for the Ph.D. degree. Supported by
  grants from the Technion-Israel Institute of Technology}
\begin{document}
\maketitle

\begin{abstract}

This paper is continuation of the paper "Primitive roots in
quadratic field". We consider an analogue of Artin's primitive
root conjecture for algebraic numbers which is not a unit in real
quadratic fields. Given such an algebraic number, for a rational
prime $p$ which is inert in the field the maximal order of the
unit modulo $p$ is $p^2-1$. An extension of Artin's conjecture is
that there are infinitely many such inert primes for which this
order is maximal. we show that for any choice of 85 algebraic
numbers satisfying a certain simple restriction, there is at least
one of the algebraic numbers which satisfies the above version of
Artin's conjecture.

\end{abstract}

\newcommand{\ord}{\operatorname{ord}}

\section{Introduction}

This paper is continuation of the paper "Primitive roots in
quadratic field" (\cite{J}). In this paper we considered an
analogue of Artin's primitive root conjecture for units in real
quadratic fields. Given such a nontrivial unit, for any rational
prime $p$ which is inert, the maximal order of the unit modulo $p$
is $p+1$. An extension of Artin's conjecture is that there are
infinitely many such inert primes for which this order is maximal.
This is known at present only under the Generalized Riemann
Hypothesis. Unconditionally, we showed that for any choice of 7
units in different real quadratic fields satisfying a certain
simple restriction, there is at least one of the units which
satisfies the above version of Artin's conjecture. In this paper
we want to extend this result for any algebraic integer modulo
inert prime $p$. we will prove\\

\begin{thm}\label{theorem 1.2}

Let $K=\mathbb{Q}(\sqrt{\Delta})$ be a quadratic field and
$\{\alpha_i\}_{i=1}^{85}$ be a set of 85 integers of K such that \\

\begin{enumerate}

\item The norms $N(\alpha_i)=\alpha_i\sigma(\alpha_i)$, of the
$\alpha_i's$, are multiplicatively  independent. \\

\item $5N(\alpha_i)\Delta, \ N(\alpha_i)$ are not perfect squares. \\

\item $ M(\alpha_i)=\sigma(\alpha_i)/\alpha_i$ are multiplicatively  independent. \\

\end{enumerate}

 Then at least one of the 85 integers has at least order $\frac{p^2-1}{24} \ \ \m \ p$
 for infinitely many inert primes $p$ in K.

\end{thm}

Note that in the case of split primes, Narkiewicz
(\cite{N1}) proved a much stronger result.\\

Since in our case the order is $p^2-1$, which is not ``linear",
some of their divisors are too big and we cannot use the method of
\cite{HB}. But since $p^2-1=(p-1)(p+1)$ can be factored into two
linear factors, with the following remark we can still use the
method of \cite{HB}

\begin{Remark} \label{Remark 1.2}

 Consider an algebraic number $\alpha$ in
 $K=\mathbb{Q}(\sqrt{\Delta})$. Let $p \nmid \alpha$ be an inert prime  in K. Since\\

 (1) $M(\alpha) \equiv \alpha^{p-1} \ (\m \ (p))$,\\

 (2) $N(\alpha)  \equiv \alpha^{p+1} \ (\m \ (p))$,\\

 we have: \\

 $\ord(M(\alpha))\mid \ord(\alpha) \ (\m \ (p))$ and $\ord(N(\alpha))\mid
 \ord(\alpha) \ (\m \ (p))$ \\

In addition: \\

 $\ord(M(\alpha))\mid p+1$ and $\ord(N(\alpha))\mid
 p-1$

 But, \\

 $(\frac{p-1}{2}, p+1)=1$ or $(p-1, \frac{p+1}{2})=1.$ \\

 So, \\

$\ord(M(\alpha))\ord(N(\alpha)) \mid 2\ord(\alpha) \ \
(\m \ (p))$ \\

\end{Remark}

Let $e_1$ and $e_2$ be some integers. If we prove that $M(\alpha)$
and $N(\alpha)$ have simultaneously at least orders
$\frac{p+1}{e_1}$ and $\frac{p-1}{e_2}$, respectively, then we
will obtain that $\alpha$ has at least order
$\frac{p^2-1}{2e_1e_2}$. This way we reduce the problem to
a ``linear" problem.\\

 \subsection{Notation and Preliminaries}

 Let $\pi(y;m,s)$ the number of primes $p \leq y$ such that
 $p\equiv s \ (\m \ m)$ where $m$ and $s$ are some integers, and
$$
E(y;m,s):= \pi(y;m,s) -\frac{Li{(y)}}{\varphi{(m)}}
$$
where $Li(y)=
 \int_2^y\frac{dt}{\log t}$. Also set
$$
E(x;m):= \max\limits_{1\leq y \leq x} \
\max\limits_{(s,m)=1}|E(y;m,s)| \;.
$$

 Denote $\mathcal{A}=\{p^2-1|p\leq x,p\equiv u \
 (\m \ {v})\}$ where $u,\ v$ are some given integers such that $(u,v)=1$, and take $X=\frac{Li(x)}{\varphi(v)}$.\\

For a square-free integer $d$, $(d,v)=1$, denote
\begin{equation*}
\begin{split}
|\mathcal{A}_d| &:=
|\{a\in\mathcal{A}:  a\equiv 0 \ \m \ d\}|\\
&=|\{p^2-1: p\leq x,\ p \equiv u \m  v,\ p^2-1 \equiv 0 \ \m \
d\}|\\&=\sum_{\substack{{m=1} \\ m^2-1 \equiv 0\m d}}^d | \{p|p
\leq x, p \equiv u \m  v,\ p \equiv m \ \m \ d\}|\\&=
\sum\limits_{\substack{{m=1} \\ m^2-1 \equiv 0 \ \m \ d \\
(m,d)=1}}^d | \{p|p \leq x, p \equiv u \ \m \ v,\ p \equiv m \ \m
\ d\}|
\end{split}
\end{equation*}

By the Chinese Remainder Theorem, for each $m$ there exists an
integer $l_m$ such that
$$
|\mathcal{A}_d| =\sum_{\substack{{m=1} \\ m^2-1 \equiv 0 \ \m \ d \\
(m,d)=1}}^d | \{p|p\leq x,\ p \equiv l_m \ (\m \ {dv})\}|;
$$

Since $|\{p|p\leq x,\ p \equiv l_m \ (\m \ {dv})\}|$ is
asymptotically independent of $m$,
 there exists some integer $l$ such that\\

$$
|\mathcal{A}_d|=\pi(x;dv,l)\sum_{\substack{m=1 \\ m^2-1 \equiv 0 \ \m \ d \\
(m,d)=1}}^d {1} = \pi(x;dv,l)\rho(d)
$$
 where $\rho(d)= \sum\limits_{\substack{m=1 \\
m^2-1 \equiv 0\ \m \ d \\ (m,d)=1}}^d {1}$.\\ \\

We note that $\rho(q)=2$ for any prime $q$, and hence for any
square-free $d$, $\rho(d)=2^{\nu(d)}$ where $\nu(d)$ denotes the
number of prime divisors of $d$. \\

 By the definition of $E(x;dv,l)$,

\begin{equation*}
\begin{split}
|\mathcal{A}_d|= \frac{\rho(d)}{\varphi(d)}\frac{Lix}{\varphi(v)}
+ \rho(d)E(x;dv,l))=\frac{2^{\nu(d)}}{\varphi(d)}X+
2^{\nu(d)}E(x;dv,l).
\end{split}
\end{equation*}

For any prime $q$ define $\omega(q):=\frac{2q}{\varphi(q)}$,
$\omega(d)=\prod\limits_{q|d}\omega(q)=\frac{2^{\nu(d)}d}{\varphi(d)}$
and
$$R_d:=|\mathcal{A}_d|-\frac{\omega(d)}{d}X = 2^{\nu(d)}E(x;dv,l)$$
Finally, we define the M$\ddot{o}$bius function, $\mu(1)=1$ and
for a
square-free $d=p_1 \cdot\cdot\cdot p_k$, $\mu(d)=(-1)^k$.\\

Now we want to prove two lemmas.\\

\begin{Lemma}\label{Lemma 8.1}
For any prime $q>3$ which is relatively prime to $v$ we have:
\begin{equation}\label{equation 8.1}
 0 \leq \frac{2}{q-1} \leq \frac{1}{2}.
 \end{equation}

\begin{equation}\label{equation 8.2}
 \sum\limits_{w \leq q < z}
\frac{2}{q-1}\log{q}-2\log\frac{z}{w} = O(1) \  \  \ (2 \leq w
\leq z)
\end{equation}
where $O$ does not depend on $z$ or $w$.

 \begin{equation}\label{equation 8.3}
\prod\limits_{\substack{2<q<z \\ q\nmid v}}
(1-\frac{2}{q-1})\gg\frac{1}{\log^2z}.
\end{equation}
 \end{Lemma}

\begin{proof}
Since $q>3$, it is clear that \eqref{equation 8.1} holds.

As for the second equality,
 $ \sum\limits_{w\leq q<z}\frac{2}{q-1}\log{q}=2\sum\limits_{w\leq q
<z}\frac{\log q}{q}\frac{q}{q-1}=2\sum\limits_{w\leq q
<z}\frac{\log q}{q}(1+\frac{1}{q-1})=2\sum\limits_{w\leq q
<z}\frac{\log q}{q}+2\sum\limits_{w\leq q <z}\frac{\log
q}{q(q-1)}=2\log\frac{z}{w}+O(1) \ \ \ \ \
(\sum\limits_{p<x}\frac{\log p}{p}= \log x+O(1))$.

Hence we get \eqref{equation 8.2}. Finally,
\begin{equation*}
\begin{split}
\prod\limits_{\substack{2<q<z \\q \nmid v}}{(1-\frac{2}{q-1})} \gg
\prod\limits_{2<q<z}{(1-\frac{2}{q-1})} \\
&= \exp({\log{\prod\limits_{2<q<z}{(1-\frac{2}{q-1})}}}) \\
&=
\exp({\sum\limits_{2<q<z}{\log{(1-\frac{2}{q-1})}}}) \\
&\gg \exp({\sum\limits_{2<q<z}(-\frac{2}{q-1}-\frac{4}{(q-1)^2})})
\end{split}
\end{equation*}

Since
$$
\frac{2}{q-1}=\frac{2}{q}+\frac{2}{q(q-1)} \leq \frac{2}{q} +
\frac{2}{(q-1)^2}
$$
and $\sum\limits_{2<q<z}\frac{2}{(q-1)^2}$ converges, we get
$$
\prod\limits_{2<q<z}(1-\frac{2}{q-1})\gg
\exp({-\sum\limits_{2<q<z}\frac{2}{q}})
$$
Since
$$\sum\limits_{2<q<z}{\frac{2}{q}}\sim 2\log\log z$$
we have
$$ \exp({-{\sum\limits_{2<q<z}\frac{2}{q}}})\gg
 \exp({-2\log\log z})=\frac{1}{\log ^2z}$$
 \end{proof}

\begin{Lemma}\label{Lemma 8.2}
For any square-free natural number $d$, $(d,v)=1$, and a real
number $A>0$, there exist constants $c_2(\geq 1)$ and $c_3(\geq
1)$ such that
 \begin{equation}\label{equation 8.4}
 \sum\limits_{d<\frac{X^\frac{1}{2}}{(\log x)^{c_2}}}
  \mu^2(d)3^{\nu(d)}|R_d| \leq c_3 \frac{X}{\log^AX}, \ \ \ (X
  \geq2)
\end{equation}
\end{Lemma}

\begin{proof}
Denote by $S_R$ the term which we need to estimate.
$$
S_R=\sum\limits_{d<\frac{X^\frac{1}{2}}{(\log
    x)^{c_2}}}\mu^2(d)3^{\nu(d)}|R_d|
$$
By the definitions of $R_d$ and $E(x;dv)$
$$
S_R\leq \sum\limits_{d<\frac{X^\frac{1}{2}}{(\log x)^{c_2}}}
\mu^2(d)6^{\nu(d)}|E(x;dv)|.
$$
Since $E(x;dv)\ll \frac{x}{dv}$ if $d\leq\frac{x}{v}$, we get that
$$
S_R \ll x^\frac{1}{2}\sum\limits_{{d<\frac{X^\frac{1}{2}}{(\log
x)^{c_2}}}}
\frac{\mu^2(d)6^{\nu(d)}}{d^\frac{1}{2}}|E(x;dv)|^{\frac{1}{2}}.
$$
 By Cauchy's Inequality,
$$
S_R\ll
x^\frac{1}{2}(\sum\limits_{d<X^\frac{1}{2}}\frac{\mu^2(d)6^{2\nu(d)}}{d})^\frac{1}{2}
(\sum\limits_{dv<\frac{vX^\frac{1}{2}}{(\log
    x)^{c_2}}}|E(x;dv)|)^\frac{1}{2}.
$$

For sufficiently large $x$ we obtain
$$
S_R\ll
x^\frac{1}{2}(\sum\limits_{d<x^\frac{1}{2}}\frac{\mu^2(d)6^{2\nu(d)}}{d})^\frac{1}{2}
(\sum\limits_{dv<\frac{x^\frac{1}{2}}{(\log
    x)^{c_2}}}|E(x;dv)|)^\frac{1}{2}.
$$

 With Bombieri-Vinogradov Theorem \cite{BV}
 (given any positive constant $e_1$, there exists a positive
 constant $e_2$ such that $\sum\limits_{d<\frac{x^\frac{1}{2}}{\log^{e_2}x}}E(x;d)=O(\frac{x}{\log^{e_1}x})$)
for the last sum and the inequality
$\sum\limits_{d<w}\frac{\mu^2(d)36^{\nu(d)}}{d}\leq (\log w
+1)^{36}$, \cite{HR}, p.115, equation (6.7)) we find that for
given
 constant $B$ there exists $c_2$ such that
$$S_R\ll \frac{x}{\log^B{x}}.$$
So, for given $A$ there exists $c_2$ such that
$$S_R \ll \frac{X}{\log^A{X}}$$
where $\ll$ depends on $v$ and $c_2$.
\end{proof}

\subsection{Proof of Theorem~\ref{theorem 1.2} -- the sieve
part} In this section we use the Selberg lower bound sieve and
show that there is some small real number $\delta_1$ and some
constant $c(\delta_1)>0$ (which depends on $\delta_1$) such that
for at least $c(\delta_1) \frac{x}{\log^3x}$ primes $p\leq x, \
p\equiv u \ (\m \ v)$, if $q|p^2-1$ then either $q >
x^{1/8+\delta_1}$ or $q|v$.

Now, define $S(\mathcal{A},z,v)=
|\{a|a\in\mathcal{A},(a,\prod\limits_{\substack{p<z \\
p \nmid v}}p)=1\}|$ and define a function $g$ by $g(t_0)=1,
t_0=4.42$ and $g(t)<1$ for $t>t_0$. Then (see \cite[Theorem 7.4,
page 219]{HR}):

\begin{Lemma}\label{Lemma 8.7} We have
\begin{equation}
 S(\mathcal{A},z,v)\geq X\prod\limits_{\substack{q<z \\ q\nmid v}} (1-\frac{\omega(q)}{q})
 \{1-g(\frac{\log X}{2\log z})+O(\frac{(loglog3X)^8 }{logX  })\}
\end{equation}
where the O-term does not depend on $X$ or on $z$.
\end{Lemma}

By Lemmas \ref{Lemma 8.1} and \ref{Lemma 8.2}, \eqref{equation
8.1}, \eqref{equation 8.2} and \eqref{equation 8.4} hold. Hence we
can use Lemma~\ref{Lemma 8.1} with $z=X^{\frac{1}{8}+\delta_0}$
$$
S(\mathcal{A},X^{\frac{1}{8}+\delta_0},v)\geq
X\prod\limits_{\substack{q<X^{\frac{1}{8}+\delta_0}\\ q\nmid v}}
(1-\frac{2}{q-1})\{1-g(\frac{1}{2}\frac{\log X}{\log
X^{\frac{1}{8}+\delta_0}}) +O(\frac{(loglog3X)^8}{logX  })\}.
$$
By Lemma~\ref{Lemma 8.1} \eqref{equation 8.3} we have for
$\delta_0$ sufficiently small
$$
S(\mathcal{A},X^{\frac{1}{8}+\delta_0},v) \gg \frac{X}{\log
  ^2X} \gg \frac{x}{\log
  ^3x}
$$

Thus for such $\delta_0$  we obtain that there is a constant
$c(\delta_0)>0$ (which depends on $\delta_0$) such that for at
least $c(\delta_0) \frac{x}{\log^3x}$ primes $p\leq x, \ p\equiv u
\ (\m \ v)$ if $q|p^2-1$ then either $q > x^{1/8+\delta_0}$ or
$q|v$. Hence we obtain that for all $0<\delta_1<\delta_0$ there is
a constant $c(\delta_1)>0$ (which depends on $\delta_1$) such that
for at least $c(\delta_1) \frac{x}{\log^3x}$ primes $p\leq x, \
p\equiv u \ (\m \ v)$, if $q|p^2-1$ then either $q >
x^{1/8+\delta_1}$ or $q|v$.

 \subsection{Proof of Theorem~\ref{theorem 1.2} -- The algebraic part}

\subsubsection{Construction of the arithmetic sequence}

Let $K=Q(\sqrt{\Delta})$ be any quadratic field, $\mathcal{O}$ the
integers ring of $K$, $ \alpha \in \mathcal{O}$ any algebraic
integer and $a=N(\alpha)$. In this section we want to construct
integers $u$ and $v$, $(u,v)=1$ such that for all primes $p$ such
that $p\equiv u \ (\m \ v)$, the discriminant $\Delta$ of
$\mathbb{Q}(\sqrt{\Delta})$ and $a$
 satisfy
$$
(\frac{\Delta}{p})=(\frac{a}{p})=-1.
$$

This means that $p$ is inert and $a$ is not a quadratic residue $
(\m \ p) $. In addition we want to obtain by the construction that
$(\frac{p^2-1}{24}, v)=1$  (since after sieving the small factors
of $\frac{p^2-1}{24}$ we may be left with small factors which
divide $v$, see previous section).

In order to fulfill these demands, we will first show that there
exist infinitely many primes $p$ satisfying the following
simultaneous conditions
\begin{equation}\label{equation 4.1}
 \textstyle(\frac{-1}{p})=(\frac{5}{p})=(\frac{a}{p})=(\frac{\Delta}{p})=-1.
\end{equation}

This condition is equivalent to the condition:
$$
B(p)=(1-(\frac{-1}{p}))(1-(\frac{5}{p}))(1-(\frac{a}{p}))(1-(\frac{\Delta}{p}))\neq
0
$$

Since the Legendre symbol is a multiplicative function, we obtain,

$$
(1-(\frac{-1}{p}))(1-(\frac{a}{p})-(\frac{\Delta}{p})+(\frac{a\Delta}{p})
-(\frac{5}{p})+(\frac{5a}{p})+(\frac{5\Delta}{p})-(\frac{5a\Delta}{p})).
$$

Let $S$ be the set of all integers of the form $n=(-1)^{b_0}
5^{b_1}a^{b_2}\Delta^{b_3}, \ \ b_i \in \{0,1\}$. Then
\begin{equation}
\textstyle\sum\limits_{p \leq Z}B(p)=\sum\limits_{n \in
S}(-1)^{b_0+b_1+b_2+b_3}\sum\limits_{p\leq Z}(\frac{n}{p}), \ \ \
b_i \in \{0,1\}.
\end{equation}

 By the assumption in the theorem each $n
\in S$ is not a perfect square when $\sum\limits_{i=0}^3b_i$ is
odd.

This assumption, together with the fact that for $n$ not a perfect
square (by the reciprocity law for Legendre symbol)
$$\textstyle\sum\limits_{p\leq Z}(\frac{n}{p})=o(\pi(Z)) \ \ as \
Z \rightarrow \infty$$ implies that $\sum\limits_{p \leq Z}B(p)$
is asymptotic to at least $\pi(Z)$ (since all the negative
summands contribute $o(\pi(Z))$ and at least the natural number 1
contributes $\pi(Z))$. This shows that the simultaneous conditions
have infinitely many solutions $p$.

We fix some particular $p_0$ satisfying the condition
\eqref{equation 4.1} and for each odd prime $l\neq 3$, such that
$l|24a\Delta$, we define $u_l=p_0$ if $l \nmid p_0^2-1$, and
$u_l=9p_0$ otherwise.
\begin{cl}
$l \nmid u_l^2-1$.
\end{cl}
 \begin{proof}If $u_l=p_0$ then by the assumption $l\nmid p_0^2-1$,
 so $l\nmid u_l^2-1$. If $u_l=9p_0$, assume, by
reductio ad absurdum, that $l|u_l^2-1$. Hence $l\mid 81p_0^2-1$.
Since, $l\mid p_0^2-1$, we obtain that $l\mid 80p_0^2$. On the
other hand, by our condition, $(\frac{5}{p_0})=-1$ so
$(\frac{p_0}{5})=-1 \ \ (p_0 \equiv 1 \ (\m \ 4))$. Hence $p_0
\equiv 2$ or $3 \ (\m \ 5)$. Since $l \mid p_0^2 -1$ and $p_0
\equiv 2$ or \ $3 \ (\m \ 5)$ we conclude that $l \nmid 5$. Using
the assumption that $l \mid p_0^2 -1$ we deduce that $ l\neq p_0 $
(if $l=p_0$ then $l \nmid p_0^2-1$). Hence ($l\neq 2,3$) $l\nmid
80p_0^2$, which is a contradiction.
\end{proof}

In addition let $u_2=p_0$ if $8 \mid {p_0}^2-1$ and $u_2=p_0-8$ if
$16 \mid {p_0}^2-1$. Likewise we take $u_3=p_0$ if $3 \mid
{p_0}^2-1$ and $u_3=p_0-3$ if $9 \mid {p_0}^2-1$. \\

 Let $v=24a\Delta$ and $u$ be the common solution of $u
\equiv u_2 \ (\m \ 16)$,  $u \equiv u_3 \ (\m \ 9)$ and all the
congruences $u\equiv u_l \ (\m \ l)$. Such  a solution exists, by
the Chinese Remainder Theorem.\\

Since $l\nmid u^2-1$ for every odd prime $l\neq 2,3, \ l\mid v$,
and by the construction $(\frac{u^2-1}{24},6)=1$ we conclude that
$(\frac{u^2-1}{24},v)=1$. \\

 Finally, if $p\equiv u \ (\m \ v)$, then $p\equiv p_0 \ (\m \ 24)$ and $p\equiv p_0$ or $4p_0
 \ (\m \ l)$ for all odd primes $l|v$. So,
$(\frac{\Delta}{p})=(\frac{\Delta}{p_0})=-1$, and similarly for
$a$. This completes the construction of $u$ and $v$. \\

Note that by the construction of the integers $u$ and $v$ we have
that $(u,v)=1$. (take $l$ an odd prime number, $l\mid v=24a\Delta$
and assume that $l\mid u$. Since $u \equiv u_l \ (\m \ l), \ \ l
\mid u_l$. Hence $l \mid p_0$ or $9p_0$ (in this case $l \neq
2,3)$. In other words $l=p_0$. But $p_0 \nmid 24a\Delta$ ($p_0$
fulfills the simultaneous condition \eqref{equation 4.1}) and $l
 \mid 24a\Delta$).

\subsubsection{The last step of the proof}

As we saw at the previous subsections, for at least $c(\delta_1)
\frac{x}{\log^3x}$ primes $p\leq x, \ p\equiv u \ (\m \ v)$, if
$q|p^2-1$ then $q > x^{1/8+\delta_1}$ or $q|v$. Since
$(\frac{p^2-1}{24},v)=1$, if $q|\frac{p^2-1}{24}$
then $q > x^{1/8+\delta_1}$.\\

Since by the construction of $u$ and $v$, $p\equiv 1 \ (\m \ 4)$
and $(\frac{p^2-1}{24},v)=1$, we have that $ \frac{p-1}{4}$ and
$\frac{p+1}{2}$ are odd. In addition if $p\equiv 1 \ (\m \ 3)$
then $(\frac{p-1}{12},v)=1$ and if
$p\equiv -1 \ (\m \ 3)$ then $(\frac{p+1}{6},v)=1$ \\

If we conclude the result about $p-1$ and $p+1$ we have for at
least $c_1(\delta_1) \frac{x}{\log^3x}$ primes $p\leq x, \ p\equiv
u \ (\m \ v)$, if $q \mid \frac{p-1}{d_-}$ or $q \mid
\frac{p+1}{d_+}$ then $q
> x^{1/8+\delta_1}$, where $d_-=4$ or 12 and $d_+=6$ or 2,
respectively.\\

For the last step of the proof we need to use a version of Lemma 4
from Narkiewicz  \cite{N1}, which generalized Lemma 2 in
\cite{GM}.
 \begin{Lemma}\label{Lemma 4.2}
If $a_1,\ldots a_k$ are multiplicatively  independent algebraic
numbers of an algebraic number-field $K$, $G$ the subgroup of
$K^\star$
 generated by $a_1,\ldots a_k$, and for any prime ideal
 $\mathbf{P}$ not dividing $a_1,\cdots a_k$ we denote by
 $G_{\mathbf{P}}$ the reduction of $G \ (\m \mathbf{P})$, then for
 all positive $y$ one can have $|G_{\mathbf{P}}| < y $ for at most
 $O(y^{1+\frac{1}{k}})$ prime ideals $\mathbf{P}$, with the implied
 constant being dependent on the $a_i$'s and $K$.
 \end{Lemma}

 \begin{proof}
 According to \cite{N1},\\

 For any real number $T$ denote by $M=M(T)$ the set of all
 k-element sequences $(r_1,...,r_k)$ of non-negative integers
 satisfying

\begin{displaymath}
 |r_1|+|r_2|+...+|r_k|\leq T
\end{displaymath}

 It is easy to see that for $T$ tending to infinity
 $|M(T)|=(c+o(1))T^k$ with suitable positive constant $c=c_k$. If
 now $P$ is a prime ideal for which $ |G_P| < y $, then select $T$ to
 be the smallest rational integer with $cT^k>2y$ and let $a_i=\frac{b_i}{c_i}$ for $i=1,...,k$. There exist two
 distinct sequences $Z=(z_i), \ W=(w_i)$ in $ M(T)$ for which

\begin{displaymath}
 P|\frac{b_1^{z_1} \cdot\cdot\cdot b_k^{z_k}}{c_1^{z_1} \cdot\cdot\cdot c_k^{z_k}}-
 \frac{b_1^{w_1} \cdot\cdot\cdot b_k^{w_k}}{b_1^{w_1} \cdot\cdot\cdot
 c_k^{w_k}}
 \end{displaymath}

  Hence for sufficient large $P$ (since the $a_i's$ are
 multiplicatively independent).

\begin{displaymath}
 P|\frac{b_1^{z_1} \cdot\cdot\cdot b_k^{z_k} - b_1^{w_1} \cdot\cdot\cdot
 b_k^{w_k}}{c_1^{[z_1,w_1]}\cdot\cdot\cdot c_k^{[z_k,w_k]}} =D \neq 0
\end{displaymath}

 Thus for sufficient large $P$

\begin{displaymath}
\nu_P(\Pi a_i^{z_i-w_i}-1) \geq 0
\end{displaymath}

where $\nu_P$ denotes the $P$-adic valuation and it follows that
for fixed $z_1-w_1,...,z_k-w_k$ we obtain
 $\ll \log{(max_j|\bar{a_j}|^{2T})} \ll T $ possibilities for $P$. Finally we
 obtain\\

$|\{ P| \ |G_P| <y \}| \ll T^{1+k} \ll y^{1+\frac{1}{k}}$

 \end{proof}

Look at the $p-1$ case (the case of $p+1$ is similar). We have for
at least $c_1(\delta_1) \frac{x}{\log^3x}$ primes $p\leq x, \
p\equiv u \ (\m \ v)$  such that if $q \mid \frac{p-1}{d_-}$ then
$q
> x^{1/8+\delta_1}$ where $d_-=4$ or 12. Let
$p-1=d_-q_1(p)q_2(p) \cdot\cdot\cdot q_m(p), \
q_m(p)>q_{m-1}(p)>...>q_1(p), \ m\leq 7$, and let $a$ be some
integer where $\bar{a}$ its image in $\mathbb{F}^*_p$.\\

 Denote by $S_n$ the set $S_n=\{a_1,...,a_n\}$ where $a_1,...,a_n$ are multiplicatively  independent
integers and take seven integers $a_{i_1},...,a_{i_7}$ from $S_n$
and assume that at least one prime, say $q_1(p)$, which is greater
than $x^{1/8+\delta_1}$, divides
$[\mathbb{F}^*_p:\langle {\bar{a}}_{i_k} \rangle ]$ for $ k=1,...,7$.\\

 Then $|\langle \bar{a}_{i_k} \rangle| = \frac{p-1}{[\mathbb{F}^*_p:\langle
\bar{a}_{i_k} \rangle]}\ll p^{7/8-\delta_1} \ll x^{7/8-\delta_1},\
\ k=1,...,7$. Since $q_1(p)$ divides $|\mathbb{F}^*_p:\langle
{\bar{a}}_{i_k} \rangle|$ for $k=1,...,7$  and $\mathbb{F}^*_p$ is
a cyclic group, $|\langle a_{i_1},...,a_{i_7} \rangle| \ \m \ p \
\ll x^{7/8-\delta_1}$. By Lemma ~\ref{Lemma 4.2} it occurs in at
most $O((x^{7/8-\delta_1})^{8/7})=O(x^{1-8/7\delta_1})$ primes
$p\leq x$ which is a negligible number relatively to
$c_1(\delta_1)
\frac{x}{\log^3x}$ for sufficiently small $\delta_1$.\\

So, for at most six integers from $a_{i_1},...,a_{i_7}, \ q_1(p)$
divides $[\mathbb{F}^*_p:\langle \bar{a_{i_k}} \rangle ]$ for
$k=1,...,7$. In other words for at most six from
$a_{i_1},...,a_{i_7} \ q_1(p)$ does not divide $|\langle
\bar{a}_{i_k} \rangle |$. Hence for at least one integer, say
$a_n$, $q_1(p)$
divide $|\langle \bar{a}_n \rangle |$.\\

Denote by $S_{n-1}$ the set $\{a_1,...,a_{n-1}\}$. By repeating
the former process for $S_{n-1}$ we obtain that for at least
one integer, say $a_{n-1}$, $q_1(p)$ divide $|\langle \bar{a}_{n-1} \rangle|$.\\

We continue this process till we obtain the set
$T_1=\{a_7,...,a_n\}$ where $a_7,...,a_n$ are multiplicatively
independent integers such that $q_1(p)$ divides $|\langle
\bar{a_{t}} \rangle|$ for $t=7,...,n$.\\

By repeating this process for $q_2(p)$ we obtain for the set
$T_2=\{a_{13},...,a_n\}$ that $q_2(p)$ divides $|\langle
\bar{a_{t}} \rangle|$ for $t=13,...,n$.\\

Again, by repeating this process for $q_m(p)$ we obtain for the
set $T_m=\{a_{6m+1},...,a_n\}$ that $q_m(p)$ divides $|\langle
\bar{a_{t}} \rangle$ for $t=6m+1,...,n$.\\

Since the maximum value of $m$ is 7, if we take $n=6m+1=43$
multiplicatively  independent integers we obtain that one of them
have at least the order $ \frac{p-1}{d_-}$ \\

Let us look on the algebraic number $ M(\alpha)$. By the same
method we obtain that one of 43 $M(\alpha)$ have at least the
order $\frac{p+1}{d_+}$. one of $42+43=85$ has, at least, the
order $\frac{p^2-1}{24}$.
 This completes the proof of Theorem ~\ref{theorem 1.2}\\

{\bf Acknowledgement} I wish to gratefully acknowledge my Ph.D
supervisors Prof. Zee`v Rudnick and Prof. Jack Sonn for the
helpful suggestion and the fruitful ideas, especially for Remark
~\ref{Remark 1.2}.

\end{document}